\documentclass[reqno,a4paper]{amsart} 
\usepackage{amsmath,amssymb,amsxtra,amscd}

\numberwithin{equation}{section} 
\theoremstyle{plain} 
\newtheorem{thm}{Theorem}[section] 
\newtheorem{cor}[thm]{Corollary} 
\newtheorem{lem}[thm]{Lemma} 
\newtheorem{Def}[thm]{Definition}

\theoremstyle{remark} 
\newtheorem{rem}[thm]{Remark} 

\newcommand\gog{\mathfrak{g}}

\begin{document}

\title{The inverse of  a real generalized Cartan matrix}

\author{Hechun Zhang} 
\address{Department of Mathematical Sciences,
 Tsinghua University, Beijing, 100084, P. R. China, and 
KdV Institute\\ Plantage Muidergracht 24\\ 1018 TV, Amsterdam,
  Netherlands\\Department of Mathematical Sciences, } \email{hzhang@math.tsinghua.edu.cn,}

\footnotetext[1]{The  author
is partially supported by NSF of China and NWO of the Netherlands.}

\dedicatory{\sc{\bf Dedicated to Professor T. H. Koornwinder on the occasion of his 60th birthday}}

\begin{abstract}
The purpose of this note is to 
give explicit criteria to determine whether a real generalized Cartan matrix  is of finite type, affine type or of hyperbolic type by considering the principal minors and the inverse of the matrix. 
In particular, it will be shown that a real generalized Cartan matrix is of finite type if and only if it is invertible and the inverse is a positive matrix. A real generalized Cartan matrix is of hyperbolic type if and only if it is invertible and the inverse is non-positive. 
\end{abstract}

\maketitle



\section{introduction}

\vskip 1cm 

For any complex $n\times n$ matrix $A$, one can associate with it a complex Lie algebra $\gog (A)$, see \cite{kac} or \cite{wan}. The properties of $A$ determine the structure of the Lie algebra $\gog (A)$. In \cite{kac}, some results are given for generalized Cartan matrices which are integral matrices. In particular, it was proved in \cite{kac} that any generalized Cartan matrix of finite type or affine type is symmetrizable. However, for any  real generalized Cartan matrix, this is not true any more, see the example  in the next section. Hence, we need to develop a new approach to deal with real generalized  Cartan matrices. In this paper, we generalize some results in \cite{kac} to arbitrary real generalized Cartan matrices. Moreover, we give some  
give explicit criteria to judge whether a real generalized Cartan matrix  is of finite type or of hyperbolic type. In particular, it will be proved that a real generalized Cartan matrix is of finite type if and only if it is invertible and the inverse is a positive matrix. A generalized Cartan matrix is of hyperbolic type if and only if it is invertible and the inverse is non-positive.
Furthermore, we obtain a similar result for real generalized Cartan matrices of affine type. 
Some of the results are even new for generalized Cartan matrices.

This work was done during the author's visit at the University of Amsterdam. The author would like to thank Professor T. H. Koornwinder for his warm hospitality and to Netherlands Organization for Scientific Research (NWO) for financial support.  The author would like also to thank Professor P. Moree for encouragement and valuable discussion.

\medskip

\section{Vinberg's classification theorem}

\vskip 1cm

The main tool of the paper is Vinberg's classification theorem. A matrix $A$ is called decomposable if, after reordering the indices (i.e. a permutation of its rows and the same permutation of its columns), $A$ decomposes into a non-trivial direct sum. A matrix that is not decomposable will be called indecomposable.

Let $u=(u_1,u_2, \cdots, u_n)^T$ be a real column vector. We write $u>0$ if $u_i>0$ for all $i$, and $u\ge 0$ if $u_i\ge 0$ for all $i$.  

\begin{Def} Let $A=(a_{ij})_{n\times n}\in M_n(\mathbb R)$. The matrix $A$ is called  a real generalized Cartan matrix if it satisfies  the following properties:

\begin{enumerate}
\item $A$ is indecomposable,
\item $a_{ij}\le 0$ for $i\ne j$,
\item $a_{ij}=0$ if and only if $a_{ji}=0$.
\end{enumerate}
\end{Def}

The following lemma was proved in \cite{kac}. 

\begin{lem} Let $A$ be a real generalized cartan matrix and let $u\ge 0$  be a real column vector. Then $Au\ge 0$ implies that either $u>0$ or $u=0$.
\end{lem}

The following theorem is known as Vinberg's  classification theorem see \cite{kac}, \cite{vin}.

\begin{thm} Let $A$ be a real generalized Cartan matrix. Then one and only one of the following possibilities holds for both $A$ and $A^T$.

\begin{enumerate}
\item $detA\ne 0$; there exists $u>0$ such that $Au>0$; $Av\ge 0$ implies $v>0$ or $v=0$. In this case, $A$ is called of finite type.
\item $detA=0$; there exists $u>0$ such that $Au=0$; $Av\ge 0$ implies $Av=0$. In this case, $A$ is called of affine type.
\item There exists $u>0$ such that $Au<0$; $Av\ge 0$, $v\ge 0$ imply $v=0$. In this case, $A$ is called of indefinite type.
\end{enumerate}
\end{thm}

It is easy to see that Vinberg's classification theorem can be simplified as follows.

\begin{cor}Let $A$ be a $n\times n$ real generalized Cartan matrix. Then
\begin{enumerate}
\item $A$ is of finite type if and only if there exists $u\in\mathbb R^n$, $u>0$ such that $Au>0$.
 \item $A$ is of affine type if and only if there exists $u\in\mathbb R^n$, $u>0$ such that $Au=0$.
\item $A$ is of indefinite type if and only if there exists $u\in\mathbb R^n$, $u>0$ such that $Au<0$.
\end{enumerate}
\end{cor}

Example. Let $A=\begin{pmatrix}4&-1&-1\\-1&4&-1\\-2&-1&4\end{pmatrix}$. Then 
$$A\begin{pmatrix}1\\1\\1\end{pmatrix}=\begin{pmatrix}2\\2\\1\end{pmatrix}.$$
Hence, $A$ is of finite type. Clearly, $A$ is non-symmetrizable.

\medskip

\section{the main results}

\vskip 1cm

Denote by $I$ the $n\times n$ identity matrix.

\begin{thm} For any real generalized Cartan matrix $A$, there exists a unique real number 
$d\in\mathbb R$ such that $A+dI$ is of affine type. Furthermore, $A+bI$ is of finite type if $b>d$ and $A+bI$ is of indefinite type  if $b<d$.\end{thm}

\begin{proof} Let
$$I_{fin}(A)=\{a\in\mathbb R| A+aI \text{ is of finite type }\},$$
$$I_{aff}(A)=\{a\in\mathbb R| A+aI \text{ is of affine type }\},$$
$$I_{ind}(A)=\{a\in\mathbb R| A+aI \text{ is of indefinite type }\}.$$

By Vinberg's classification theorem, we have

$$\mathbb R=I_{fin}(A)\dot{\cup} I_{aff}(A)\dot{\cup} I_{ind}(A).$$

 If $a\in I_{fin}(A)$ and $b\ge a$, there exists $u>0$ such that $(A+aI)u>0$. Hence, $(A+bI)u>0$ and so  $b\in I_{fin}(A)$. For this given real number $a$, we can find a sufficient small $\epsilon>0$ such that $(A+(a-\epsilon)I)u>0$.  Hence, $a-\epsilon\in I_{fin}(A)$ and $I_{fin}(A)$ is an open set.
Similarly, if $a\in I_{ind}(A)$ and $b\le a$, then there exists $u>0$ such that $(A+aI)u<0$, then $(A+bI)u<0$ and so $b\in I_{ind}(A)$. For this given real number $a$, we can find a sufficient small $\epsilon>0$ such that $(A+(a+\epsilon)I)u<0$. Hence, $a+\epsilon\in I_{ind}(A)$ and $I_{ind}(A)$ is an open set.

If $a\in I_{aff}$, then $-a$ is a real eigenvalue of the matrix $A$. Hence, $I_{aff}(A)$ is a finite set. Since $I_{fin}(A)$ and $I_{ind}(A)$ are open sets, $I_{aff}(A)$ is non-empty. Assume that $a_1, a_2\in I_{aff}(A)$ and $a_1<a_2$. Then there exists $b$, $a_1<b<a_2$, such that $A+bI$ is of finite type or indefinite type. In either cases, it is a contradiction. Hence, $I_{aff}(A)$ only contains one element, say $d$. Moreover,
$$I_{fin}(A)=\{a\in\mathbb R|a>d\}$$ 
and 
$$I_{ind}(A)=\{a\in\mathbb R|a<d\}.$$
\end{proof}

The following lemma was proved in \cite{kac}.

\begin{lem} Let $A$ be a real generalized Cartan matrix. If $A$ is of finite or affine type, then any proper principal submatrix decomposes into a direct sum of matrices of finite type.\end{lem}

\begin{thm} The real generalized Cartan matrix $A$ is of finite type if and only if $detA\ne 0$ and $A^{-1}>0$.\end{thm}

\begin{proof} Let $A$ be a real generalized  Cartan matrix of finite type. 
Then $detA\ne 0$ by Vinberg's theorem. Let $A^{-1}=(u_1, u_2,\cdots, u_n)$ where $u_i$'s are the column vectors of $A^{-1}$. Then $Au_i\ge 0$ for all $i$. Again by Vinberg's theorem, we have $u_i>0$ for all $i$, i.e. $A^{-1}>0$.

Conversely, assume that $A$ is a generalized Cartan matrix such that $detA\ne 0$ and $A^{-1}>0$. Of course, $A$ is not of affine type. If $A$ is of indefinite type, then there exists $u>0$ such that $Au=w<0$. Since $A^{-1}>0$, we have $u=A^{-1}w<0$ which is a contradiction. Therefore, $A$ is of finite type. 
\end{proof}

\begin{Def} A real generalized Cartan matrix $A$ is called a matrix of hyperbolic type if it is of indefinite type and any proper principal submatrix of $A$ is a direct sum of matrices of  finite type or affine type. A hyperbolic generalized Cartan matrix $A$ is called strictly hyperbolic if any proper principal submatrix of $A$ is a direct sum of matrices of finite type. 
\end{Def}

\begin{thm}\label{hy1} Let $A$ be a real generalized Cartan matrix. Then $A$ is of hyperbolic type if and only if 
for any column vector $u\in\mathbb R^n$, $Au\le 0$ implies $u\ge 0$.
\end{thm} 

\begin{proof}

Let $A$ be a generalized Cartan matrix of hyperbolic type. Let $0\ne u\in \mathbb R^n$ such that $Au\le 0$. We claim that $u\ge 0$. If $u$ is not non-negative, then $u$ has both positive and negative entries since $u\le 0$ contradicts to that $A$ is of indefinite type. By reordering, we may assume that 
$$u_1>0,\cdots, u_{s-1}>0, u_s\le 0,\cdots, u_{n-1}\le 0, u_n<0.$$
We divide $A$ into block form

$$A=\begin{pmatrix}A_1&A_2\\A_3&A_4\end{pmatrix}$$
where $A_1$ is a $(s-1)\times (s-1)$ submatrix, $A_2$ is a $(s-1)\times (n-s+1)$ submatrix, $A_3$ is a $(n-s+1)\times (s-1)$ submatrix and $A_4$ is a  
$(n-s+1)\times (n-s+1)$ submatrix.  Let $u^{\prime}=(u_1,\cdots, u_{s-1})^T$ and let $w^{\prime}=(u_s,\cdots, u_n)^T$. Then we have

$$A_1u^{\prime}+A_2w^{\prime}\le 0.$$
Hence, $A_1u^{\prime}\le -A_2w^{\prime}\le 0$. Since $A_1$ is a direct sum of matrices of finite type or affine type, we conclude that $A_1$ is of affine type and so $s=n$ and $A_2u_n=0$. Therefore, $A_2=0$ and so $A_3=0$ which contradicts to that $A$ is indecomposable. Hence, $u\ge 0$.

If $A$ is of indefinite type but is not of hyperbolic type, then $A$ has a $(n-1)\times (n-1)$ principal submatrix $A_1$ which is of indefinite type. We may assume that 

$$A=\begin{pmatrix}A_1&-b_1\\-c_1^T&a_{nn}\end{pmatrix},$$
where $b_1\ge 0, c_1\ge 0$ are column vectors in $\mathbb R^{n-1}$. There exists $u>0$ in $\mathbb R^{n-1}$ such that $A_1u<0$. For $a\in \mathbb R_+$, set
$$w=\begin{pmatrix}au\\-1\end{pmatrix}.$$
Then we have $Aw<0$ for sufficient large $a$. Contradiction! \end{proof}

\begin{thm} Let $A$ be a  real generalized Cartan matrix. Then $A$ is of hyperbolic type if and only if $detA\ne 0$ and $A^{-1}\le 0$.\end{thm}

\begin{proof} Assume that  $A$ is of hyperbolic type. If $detA=0$, then there exist $0\ne u\in \mathbb R^n$ such that $Au=A(-u)=0$. By Theorem \ref{hy1}, $u\ge 0$ and $-u\ge 0$ which force $u=0$. Hence, $A$ is nondegenerate. Assume that 
$$A^{-1}=(u_1, u_2,\cdots, u_n),$$
where $u_1, u_2,\cdots, u_n\in\mathbb R^n$. Then $Au_i\ge 0$ for all $i$. By theorem \ref{hy1}, $u_i\le 0$. 

Conversely, assume that $detA\ne 0$ and $A^{-1}\le 0$. Let $u\in\mathbb R^n$ be such  that $Au=w\ge 0$. Then $u=A^{-1}w\le 0$. Again by Theorem \ref{hy1}, $A$ is of hyperbolic type.
\end{proof}

\begin{thm} The real generalized Cartan matrix $A$ is of finite  type if and only if all of the principal minors of $A$ are positive. The real generalized Cartan matrix $A$ is of affine type if and only if $detA=0$ and all of the principal minors with size $\le n-1$ are positive.\end{thm}

\begin{proof}
Let $B$ be an arbitrary real generalized Cartan matrix of finite type. Then for any $a\ge 0$, $B+aI$ is also of finite type. Hence, $det(B+aI)\ne 0$ which means the real eigenvalues (if there are any) of $B$ are positive. Therefore, $detB>0$. 

Now assume that $A=(a_{ij})_{n\times n}$ is a generalized Cartan matrix and that  all of the principal minors of $A$ are positive. 
We then use induction on $n$ to show that the matrix $A$ is of finite type. The case $n=1$ is trivial. If $A$ is not of finite type, then $A$ must be of indefinite type because $A$ can not be of affine type. By the induction hypothesis, any principal indecomposable submatrix is of finite type. Hence, $A$ is of strictly hyperbolic type. Therefore, $A^{-1}\le 0$. But the $(1,1)$-entry of $A^{-1}$ is a $(n-1)\times (n-1)$ principal minor of $A$ over $detA$ which means that $detA<0$. Contradiction! Hence, $A$ must be of finite type.

If $A$ is of affine type, then  any proper principal submatrix  of $A$ is a direct sum of matrices of finite type. Hence, the determinant of the proper principal matrix is positive. 

Conversely, assume that  $detA=0$ and all of the principal minors with size $\le (n-1)$ are all positive. Clearly, $A$ is not of finite type. If $A$ is of indefinite type, then $A$ must be of strictly hyperbolic type, and so $detA\ne 0$ which is a contradiction. Therefore, $A$ is of affine type.
\end{proof}

\begin{rem}If a real generalized Cartan matrix $A$ is of finite type (resp. affine type) and is symmetrizable, then it is easy to see that there exist positive real numbers $d_1, d_2, \cdots, d_n$ (see \cite{kac} 2.3) and a symmetric matrix $B$ such that 
$$A=diag(d_1, d_2, \cdots, d_n)B.$$
From the above theorem, one can see that $B$ is positive definite (resp. semi positive definite).
\end{rem}

\begin{cor}Let $A$ be a RGCM of hyperbolic type. Then $detA<0$.\end{cor}

\begin{proof} We already proved that $detA<0$ if $A$ is of strictly hyperbolic type. Assume that $A$ is of hyperbolic type. For a sufficient small $\epsilon>0$, $A+\epsilon I$ is still of indefinite type. However, any principle indecomposable  submatrix is of finite type. Hence, $det(A+\epsilon I)<0$. Note that $detA\ne 0$. Let $\epsilon$ tend 
to zero. Then $detA<0$. 
\end{proof}

\bibliographystyle{amsplain}

\end{document}